\documentclass[11pt]{article}
\usepackage{mathrsfs}
\usepackage{amssymb}
\usepackage{amsfonts}
\usepackage{color,xcolor}
\usepackage{graphicx}
\usepackage{geometry}
\usepackage{manfnt}
\usepackage[pagewise]{lineno}
\usepackage[colorlinks,citecolor=blue,urlcolor=blue]{hyperref}
\newtheorem{theorem}{Theorem}[section]
\newtheorem{corollary}{Corollary}[section]
\newtheorem{lemma}{Lemma}[section]

\newtheorem{remark}{Remark}[section]

\def\[{{\Big[}}\def\]{{\Big]}}\def\({{\Big(}}\def\){{\Big)}}

\def\cC{{\mathcal C}}

\def\mE{{\mathbb E}}
\def\mK{{\mathbb K}}
\def\mP{{\mathbb P}}
\def\mR{{\mathbb R}}

\def\={&\!\!=\!\!&}

\def\geq{\geqslant}\def\leq{\leqslant}

\begin{document}
\title{\bf Strong solutions of stochastic differential equations with square integrable drift}

\author{Rongrong Tian$^a$, Liang Ding$^b$ and Jinlong Wei$^c$\\
{\small \it $^a$College of Science, Wuhan University of Technology, Wuhan 430070, China}
\\ {\small \tt
tianrr2018@whut.edu.cn}
\\ {\small \it $^b$School of Data Science and Information Engineering, Guizhou Minzu } \\  {\small \it University, Guiyang, 550025, China}\\
{\small \tt ding2016liang@126.com}\\
{\small \it $^c$School of Statistics and Mathematics, Zhongnan University of}\\ {\small \it Economics and Law, Wuhan 430073, China}  \\
 {\small \tt weijinlong.hust@gmail.com} }
\date{}

\maketitle

\noindent{\hrulefill}
\vskip1mm\noindent{\bf Abstract} We prove the existence and uniqueness of strong solutions for stochastic differential equations in which the drift coefficient is square integrable in time variable and H\"{o}lder continuous in space variable.  Moreover, we prove that the unique strong solution has a continuous modification, which is $\beta$-H\"{o}lder continuous in space variable for every $\beta\in (0,1)$, and as an $L^2(\Omega\times (0,T))$ valued function, it is differentiable as well.


\vskip1mm\noindent{\hrulefill}

\section{Introduction}\label{sec1}\setcounter{equation}{0}
Consider the following stochastic differential equation (SDE for short) in $\mathbb{R}^d$:
\begin{eqnarray}\label{1.1}
dX_t=b(t,X_t)dt+\sigma(t,X_t)dW_t, \ t\in (0,T], \
X_0=x\in{\mathbb R}^d,
\end{eqnarray}
where $T>0$ is a given real number, $\{W_t\}_{0\leq t\leq T}=\{(W_{1,t}, \ldots, W_{d,t})\}_{0\leq t\leq T}$ is
a $d$-dimensional standard Wiener process defined on a given stochastic
basis ($\Omega, {\mathcal F},\{{\mathcal F}_{t}\}_{0\leq t\leq T},{\mathbb P}$), and the coefficients $b: [0,T]\times{\mathbb R}^d\rightarrow{\mathbb R}^d$, $\sigma: [0,T]\times{\mathbb R}^d\rightarrow{\mathbb R}^{d\times d}$
are Borel measurable. When $\sigma$ is Lipschitz continuous in $x$
uniformly in $t$ and $b$ is bounded measurable, Veretennikov \cite{Ver}
first proved the existence of a unique strong solution for SDE (\ref{1.1}). Since then,
Veretennikov's result was strengthened in different forms under the same
assumption on $b$, see \cite{Dav07,Fla11,MNP}. When
$\sigma=I_{d\times d}$ and $b$ is in the Krylov-R\"{o}ckner class:
\begin{eqnarray}\label{1.2}
b\in L^q(0,T;L^p(\mathbb{R}^d;\mathbb{R}^d)), \ \ p,q\in[2,\infty),
\end{eqnarray}
and
\begin{eqnarray}\label{1.3}
\frac{2}{q}+\frac{d}{p}<1,
\end{eqnarray}
using the Girsanov transformation and Krylov estimates, Krylov and
R\"{o}ckner \cite{KR} proved the existence and uniqueness of strong
solutions for (\ref{1.1}). Recently, following \cite{KR}, Fedrizzi
and Flandoli \cite{FF13} derived that the unique strong solution forms a $\cC^\beta$ ($\beta\in (0,1)$) stochastic homeomorphism flow. More recently, for non-constant diffusion,
if $\sigma(t,x)$ is continuous in $x$ uniformly with respect to $t$,
$\sigma\sigma^\top$ satisfies a uniformly elliptic condition and
$|\nabla\sigma|\in L^q(0,T;L^p(\mathbb{R}^d))$ with $p,q$
satisfying (\ref{1.3}), Zhang \cite{Zha11} demonstrated the existence and
local uniqueness for the stochastic homeomorphism flow. Moreover, there are
many other excellent research works devoted to studying the stochastic homeomorphism flow for SDE (\ref{1.1}) under various non-Lipschitz coefficients,
we refer to \cite{Att,Fang,Luo09,Luo15,Qiao,Ren03,Ren10}.

However, there are few investigations to deal with the case of
(\ref{1.2}) with $q\leq2$, and now the condition (\ref{1.3}) is no longer true.
To get an analogue condition of (\ref{1.3}), we assume that
\begin{eqnarray}\label{1.4}
b\in L^q(0,T;\mathcal{C}^\alpha_b(\mathbb{R}^d;\mathbb{R}^d)), \ \alpha\in (0,1).
\end{eqnarray}
When $q=\infty$ and $\sigma=I_{d\times d}$, the existence as well as
uniqueness for strong solution has been established by Flandoli, Gubinelli
and Priola \cite{FGP1} (also see \cite{FGP2} for unbounded  drift
coefficients). This result has been generalized by Wei, Duan, Gao and Lv \cite{WDGL} to $q>2/\alpha$ (also see \cite{Tian1}).  As we known, there are still few research works concerned with $q\leq 2$. This problem is the main driving source for us to work out the present paper. Under the condition (\ref{1.4}) with $q\leq 2$, we will prove the existence, uniqueness, H\"{o}lder continuity and weak
differentiability of solutions for (\ref{1.1}). Let $\mathcal{C}_{u}(\mathbb{R}^d)$ be the space consisting of all bounded uniformly continuous functions on $\mathbb{R}^d$. Then our main result is given as the following:
\begin{theorem}\label{thm1.1} Let $b\in L^1(0,T;\mathcal{C}_{u}(\mathbb{R}^d;\mathbb{R}^d))$, and
let $\sigma=(\sigma_{i,j})$ be a $d\times d$ matrix valued function such that $\sigma_{i,j}\in L^2(0,T;\mathcal{C}_{u}(\mathbb{R}^d))$ ($1\leq i,j\leq d$).

(i) Then there exists a weak solution to SDE (\ref{1.1}), i.e. there is a filtered
probability space
$(\tilde{\Omega},\tilde{\mathcal{F}},\{\tilde{\mathcal{F}}_t\}_{0\leq t\leq T},\tilde{\mathbb{P}})$,
two processes $\tilde{X}_t$ and $\tilde{W}_t$ defined for $t\in [0,T]$ on
it such that $\tilde{W}_t$ is a $d$-dimensional
$\{\tilde{\mathcal{F}}_t\}$-Wiener process and $\tilde{X}_t$ is an
$\{\tilde{\mathcal{F}}_t\}$-adapted, continuous, $d$-dimensional process
and for every $t\in [0,T]$,
\begin{eqnarray}\label{1.5}
\tilde{X}_t=x+\int_0^tb(r,\tilde{X}_r)dr+
\int_0^t\sigma(r,\tilde{X}_r)d\tilde{W}_r, \ \ \mP-a.s..
\end{eqnarray}

(ii) Let $\alpha\in (0,1)$ and let (\ref{1.4}) hold. For $1\leq i,j\leq d$, we assume that
$\sigma_{i,j}\in L^\infty(0,T;\mathcal{C}_b^\alpha(\mathbb{R}^d))$ and $\nabla\sigma_{i,j}\in L^2(0,T;L^\infty(\mathbb{R}^d;\mR^d))$. We assume further that $a=(a_{i,j})=\sigma\sigma^\top=(\sigma_{i,k}\sigma_{j,k})$ is uniformly
elliptic, i.e. for every $(t,x)\in (0,T)\times \mathbb{R}^d$ there is a constant $\Lambda>1$ such that
\begin{eqnarray}\label{1.6}
\frac{1}{\Lambda} |\xi|^2\leq \sum_{i,j=1}^da_{i,j}(t,x)\xi_i\xi_j\leq \Lambda |\xi|^2, \
\forall \ \xi=(\xi_1,\xi_2,\ldots,\xi_d)\in\mathbb{R}^d.
\end{eqnarray}
Then the pathwise uniqueness holds for (\ref{1.1}). Moreover, the stochastic field $\{X_t(x), t\in [0,T],x\in{\mathbb R}^d\}$ has a continuous modification $\tilde{X}$,  which is $\beta$-H\"{o}lder continuous in $x$ for every $\beta\in (0,1)$, and for every $p\geq 1$,
\begin{eqnarray}\label{1.7}
\mathbb{E}\Big[\sup_{0\leq t\leq T}\Big(\sup_{x\neq y}\frac{|\tilde{X}_t(x)-\tilde{X}_t(y)|}{|x-y|^\beta}\Big)^p\Big]<\infty.
\end{eqnarray}

(iii) With the same conditions of (ii). For almost surely $\omega\in\Omega$,
every $t\in [0,T]$, $x\mapsto X_t(x)$ is a homeomorphism on $\mathbb{R}^d$. Moreover,
$X_t(x)$ is differentiable in $x$ in the sense that: $\{e_i\}_{i=1}^{d}$ is
the canonical basis of ${\mathbb R}^d$, for every $x\in {\mathbb R}^d$ and $1\leq i\leq d$, the limit
\begin{eqnarray}\label{1.8}
\lim_{\delta\rightarrow 0}\frac{X_\cdot(x+\delta e_i)-X_\cdot(x)}{\delta}
\end{eqnarray}
exists in $L^2(\Omega\times (0,T))$.
\end{theorem}

In this paper, when there is no ambiguity,
we use $C$ to denote a constant whose true value may vary from line to line. We
use $\nabla$ to denote the gradient of a function with respect to the space variable.  As usual, ${\mathbb N}$ stands for
the set of all natural numbers.  a.s. is the abbreviation of almost surely.

\section{Schauder estimates for parabolic partial differential equations}\label{sec2}
\setcounter{equation}{0}
Assume $\alpha\in(0,1)$, $T>0$ and $q\in[1,2]$. Let $L^q(0,T;\mathcal{C}^{\alpha}_{b}(\mathbb{R}^d))$
denote the set consisting of all bounded Borel functions $f: [0,T]\times \mR^d\rightarrow\mR$ for which
\begin{eqnarray*}
\|f\|_{L^q(0,T;\mathcal{C}^\alpha_b(\mathbb{R}^d))}
&=&\Big\{\int_0^T\Big[ \sup_{x\in\mathbb{R}^d}|f(t,x)|+\sup_{x,y\in\mathbb{R}^d, x\neq y}\frac{|f(t,x)-f(t,y)|}{|x-y|^\alpha}\Big]^qdt\Big\}^{\frac{1}{q}}
\nonumber\\&=&:
\Big\{\int_0^T\Big[\|f(t)\|_0+[f(t)]_\alpha \Big]^qdt\Big\}^{\frac{1}{q}}<\infty.
\end{eqnarray*}
Moreover, for $n\geq1$, $f\in L^q(0,T;\mathcal{C}^{n,\alpha}_{b}(\mathbb{R}^d))$ if all spatial partial derivatives $\partial^{k}_{x_{i_1},\ldots,x_{i_k}}f$ are in $L^q(0,T;\mathcal{C}^{\alpha}_{b}(\mathbb{R}^d))$, for all orders $k=0,1,\ldots,n$ and $1\leq i_1,\ldots,i_k\leq d$. The norm in $L^q(0,T;\mathcal{C}^{n,\alpha}_{b}(\mathbb{R}^d))$ is defined by
\begin{eqnarray*}
\|f\|_{L^q(0,T;\mathcal{C}^{n,\alpha}_b(\mathbb{R}^d))}
=\Big\{\int_0^T\Big[\sum_{k=0}^n\|\partial^{k}_{x_{i_1},\ldots,x_{i_k}}f(t)\|_0+[\partial^n_{x_{i_1},\ldots,x_{i_n}}f(t)]_\alpha \Big]^qdt\Big\}^{\frac{1}{q}}.
\end{eqnarray*}
Similarly, we define the space $W^{1,q}(0,T;\mathcal{C}^{\alpha}_{b}(\mathbb{R}^d))$ as the set consisting of all bounded Borel functions $f: [0,T]\times \mR^d\rightarrow\mR$ for which $f,\partial_t f\in L^q(0,T;\mathcal{C}^{\alpha}_{b}(\mathbb{R}^d))$, the norm is defined by
\begin{eqnarray*}
\|f\|_{W^{1,q}(0,T;\mathcal{C}^\alpha_b(\mathbb{R}^d))}
=\|f\|_{L^q(0,T;\mathcal{C}^\alpha_b(\mathbb{R}^d))}+\|\partial_tf\|_{L^q(0,T;\mathcal{C}^\alpha_b(\mathbb{R}^d))}.
\end{eqnarray*}
If a function $f$ is in $L^q(0,T;\mathcal{C}^{n,\alpha}_{b}(\mathbb{R}^d))\cap W^{1,q}(0,T;\mathcal{C}^{\alpha}_{b}(\mathbb{R}^d))$, we set
\begin{eqnarray*}
\|f\|_{L^q(0,T;\mathcal{C}^{n,\alpha}_{b}(\mathbb{R}^d))\cap W^{1,q}(0,T;\mathcal{C}^\alpha_b(\mathbb{R}^d))}
=\|f\|_{L^q(0,T;\mathcal{C}^{n,\alpha}_b(\mathbb{R}^d))}+
\|f\|_{W^{1,q}(0,T;\mathcal{C}^\alpha_b(\mathbb{R}^d))}.
\end{eqnarray*}

Let $g$ and $h$ be Borel measurable. Consider the following Cauchy problem for $u:[0,T]\times\mathbb{R}^d\to\mathbb{R}$,
\begin{eqnarray}\label{2.1}
\left\{\begin{array}{ll}
\partial_{t}u(t,x)=\frac{1}{2}\Delta u(t,x)+g(t,x)\cdot \nabla u(t,x)+h(t,x), \ (t,x)\in (0,T)\times {\mathbb R}^d, \\
u(0,x)=0, \  x\in{\mathbb R}^d.  \end{array}\right.
\end{eqnarray}
$u$ is said to be a strong solution of (\ref{2.1}) if it is
in $L^q(0,T;\mathcal{C}^{2,\alpha}_b(\mathbb{R}^d))\cap W^{1,q}(0,T;\mathcal{C}^{\alpha}_b(\mathbb{R}^d))$ such that (\ref{2.1}) holds true for almost all $(t,x)\in (0,T)\times \mR^d$.

Now, we give a useful lemma, which claims that the strong solution has an equivalent form.
\begin{lemma} \label{lem2.1} Suppose
$g\in L^q(0,T;\mathcal{C}^{\alpha}_b(\mathbb{R}^d;\mathbb{R}^d))$ and $h\in L^q(0,T;\mathcal{C}^{\alpha}_b(\mathbb{R}^d))$. Let $u\in L^q(0,T;\mathcal{C}^{2,\alpha}_b(\mathbb{R}^d))\cap W^{1,q}(0,T;\mathcal{C}^{\alpha}_b(\mathbb{R}^d))$. If $u$ is a strong solution of (\ref{2.1}), then for every $t\in [0,T]$ and almost all $x\in\mathbb{R}^d$,
it fulfils the integral equation
\begin{eqnarray}\label{2.2}
u(t,x)=\int_0^tK(r,\cdot)\ast (g(t-r,\cdot)\cdot \nabla u(t-r,\cdot))(x)dr+\int_0^tK(r,\cdot)\ast h(t-r,\cdot)(x)dr,
\end{eqnarray}
where $K(r,x)=(2\pi r)^{-\frac{d}{2}}e^{-\frac{|x|^2}{2r}}, r>0, x\in\mathbb{R}^d$, and vice versa.
\end{lemma}
\smallskip
\vskip2mm\noindent
\textbf{Proof.} If (\ref{2.2}) holds, by using the properties for the heat kernel $K$, then (\ref{2.1}) is true. On the other hand, if $u$ satisfies (\ref{2.1}), then for every $\psi\in \mathcal{C}_0^\infty({\mathbb R}^d)$
\begin{eqnarray*}
\int_0^t\int_{\mathbb{R}^d}\partial_ru(r,x)\varphi(r,x)drdx&=&
\frac{1}{2}\int_{\mathbb{R}^d}\int_0^t\Delta u(r,x)\varphi(r,x)drdx+ \int_{\mathbb{R}^d}\int_0^t
g(r,x)\cdot \nabla u \varphi(r,x) drdx \nonumber\\&&+\int_{\mathbb{R}^d}\int_0^t h(r,x)\varphi(r,x)drdx, \quad \forall \ t\in [0,T],
\end{eqnarray*}
where $\varphi(r,x)=K(t-r,\cdot)\ast \psi(\cdot)(x)$.

In view of integration by parts, it follows that
\begin{eqnarray*}
\int_{\mathbb{R}^d}u(t,x)\psi(x)dx&=&
\int_{\mathbb{R}^d}\int_0^t u(r,x)[\partial_r\varphi(r,x)+\frac{1}{2}\Delta\varphi(r,x)]drdx \nonumber\\&&+ \int_{\mathbb{R}^d}\int_0^t
g(r,x)\cdot \nabla u(r,x) \varphi(r,x) drdx+\int_{\mathbb{R}^d}\int_0^t h(r,x)\varphi(r,x)drdx
\nonumber\\ &=& \int_{\mathbb{R}^d}\int_0^tK(t-r,\cdot)\ast (g(r,\cdot)\cdot \nabla u(r,\cdot))(x)dr\psi(x)dx\nonumber\\&&+\int_{\mathbb{R}^d}\int_0^tK(t-r,\cdot)\ast h(r,\cdot)(x)dr\psi(x)dx, \quad \forall \ t\in [0,T].
\end{eqnarray*}
Since $\psi$ is arbitrary, (\ref{2.2}) is true. $\Box$

In preparation to the next section, the following lemma will play an
important role.
\begin{lemma} \label{lem2.2} Let $q,\alpha,h$ and $g$ be stated in
Lemma \ref{lem2.1}. We assume further that
\begin{eqnarray}\label{2.3}
\theta:=1+\alpha-2/q>0.
\end{eqnarray}

(i) If $g=0$, (\ref{2.1}) has a unique strong solution $u$. Moreover $u\in {\mathcal C}([0,T];{\mathcal C}_b^{1,\theta}({\mathbb R}^d))$ and
\begin{eqnarray}\label{2.4}
\|u\|_{{\mathcal C}([0,T];{\mathcal C}_b^{1,\theta}({\mathbb R}^d))}
\leq C \|h\|_{L^q(0,T;\mathcal{C}^\alpha_b(\mathbb{R}^d))}.
\end{eqnarray}

(ii) For a general function $g$, we assume $q=2$ in addition,
then (\ref{2.1}) exists a unique strong solution $u$.
Moreover,
$u$ belongs to ${\mathcal C}([0,T];{\mathcal C}_b^{1,\alpha}({\mathbb R}^d))$ and
\begin{eqnarray}\label{2.5}
\|u\|_{{\mathcal C}([0,T];{\mathcal C}_b^{1,\alpha}({\mathbb R}^d))} \leq C\|h\|_{L^2(0,T;\mathcal{C}^{\alpha}_b(\mathbb{R}^d))},
\end{eqnarray}
and
\begin{eqnarray}\label{2.6}
\|u\|_{L^2(0,T;\mathcal{C}^{2,\alpha}_b(\mathbb{R}^d))\cap W^{1,2}(0,T;\mathcal{C}^{\alpha}_b(\mathbb{R}^d)) } \leq C\|h\|_{L^2(0,T;\mathcal{C}^{\alpha}_b(\mathbb{R}^d))},
\end{eqnarray}
where the constant $C$ in (\ref{2.5}) and (\ref{2.6}) depends only on $\alpha,d$ and $\|g\|_{L^2(0,T;\mathcal{C}^\alpha_b(\mathbb{R}^d;\mR^d))}$.
\end{lemma}
\vskip2mm\noindent
\textbf{Proof.} $(i)$ Let $u$ be given by (\ref{2.2}). Using Hausdorff-Young's inequality, then $u\in L^q(0,T;W^{1,\infty}(\mathbb{R}^d))$ and there is a constant $C>0$ which only depends on $\alpha,q,d$ and $T$ such that
 \begin{eqnarray}\label{2.7}
\int_0^T\|u(t)\|^q_{W^{1,\infty}(\mathbb{R}^d)}dt \leq C\int_0^T\|h(t)\|_0^qdt.
\end{eqnarray}
From (\ref{2.3}), we know $q>1$. With the help of \cite[Theorem 3.3]{Kry02}, for every $1\leq i,j\leq d$,
$\partial^2_{x_i,x_j}u\in L^q(0,T;\mathcal{C}^\alpha_b(\mathbb{R}^d))$ and there is another constant $C>0$ such that
\begin{eqnarray}\label{2.8}
\int_0^T[\partial^2_{x_i,x_j}u(t)]^q_\alpha dt \leq C\int_0^T[h(t)]^q_\alpha dt.
\end{eqnarray}
Using (\ref{2.7}) and (\ref{2.8}), we infer that $u\in L^q(0,T;\mathcal{C}^{2,\alpha}_b(\mathbb{R}^d))$.
On the other hand, $u$ satisfies (\ref{2.1}) in the distributional sense, so it also lies in $W^{1,q}(0,T;\mathcal{C}^{\alpha}_b(\mathbb{R}^d))$.
By Lemma \ref{lem2.1}, (\ref{2.1}) exists a unique strong solution. It remains to  prove $u\in {\mathcal C}([0,T];{\mathcal C}_b^{1,\theta}({\mathbb R}^d))$ and (\ref{2.4}).
Clearly, it needs to check: for every $1\leq i\leq d$,
\begin{eqnarray}\label{2.9}
\partial_{x_i}u(\cdot,\cdot)=\partial_{x_i}\Big[\int_0^{\cdot}K(r,\cdot)\ast
h(\cdot-r,\cdot)dr\Big]\in{\mathcal C}([0,T];{\mathcal C}_b^\theta({\mathbb R}^d)),
\end{eqnarray}
and the inequality (\ref{2.4}).

Let $x,y\in{\mathbb R}^d$ and let $t\in [0,T)$. For $\delta>0$ ($t+\delta<T$) and
$1\leq i\leq d$, we divide the quantity
\begin{eqnarray*}
[\partial_{x_i}u(t+\delta,x)-\partial_{x_i}u(t,x)]-
[\partial_{y_i}u(t+\delta,y)-\partial_{y_i}u(t,y)]
\end{eqnarray*}
into $\sum_{i=1}^8I_i^\delta(t)$ with
\begin{eqnarray*}
&&I_1^\delta(t)=\int_0^tdr\int_{|x-z|\leq 2|x-y|}\partial_{x_i}K(r,x-z)[h^\delta(t-r,z)-h^\delta(t-r,x)]dz,\nonumber\\
\nonumber\\
&&I_2^\delta(t)=\int_0^tdr\int_{|x-z|\leq 2|x-y|}\partial_{y_i}K(r,y-z)
[h^\delta(t-r,y)-h^\delta(t-r,z)]dz,\nonumber\\
&&I_3^\delta(t)=\int_0^tdr\int_{|x-z|> 2|x-y|}\partial_{y_i}K(r,y-z)[h^\delta(t-r,y)-h^\delta(t-r,x)]dz,
\nonumber\\
&&I_4^\delta(t)=\int_0^tdr\int_{|x-z|> 2|x-y|}[\partial_{x_i}K(r,x-z)-\partial_{y_i}K(r,y-z)][h^\delta(t-r,z)
-h^\delta(t-r,x)]dz,
\nonumber\\
&&I_5^\delta(t)=
\int_t^{t+\delta}dr\int_{|x-z|\leq 2|x-y|}\partial_{x_i}K(r,x-z)h^\delta(t-r,z,x)dz,\nonumber\\
&&I_6^\delta(t)=\int_t^{t+\delta}dr\int_{|x-z|\leq 2|x-y|}\partial_{y_i}K(r,y-z)h^\delta(t-r,y,z)dz,
\nonumber\\
&&I_7^\delta(t)=\int_t^{t+\delta}dr\int_{|x-z|> 2|x-y|}\partial_{y_i}K(r,y-z)h^\delta(t-r,y,x)dz,\nonumber\\
&&I_8^\delta(t)=\int_t^{t+\delta}dr\int_{|x-z|> 2|x-y|}[\partial_{x_i}K(r,x-z)-\partial_{y_i}K(r,y-z)] h^\delta(t-r,z,x)dz,
\end{eqnarray*}
where
$$
h^\delta(t-r,x)=h(t+\delta-r,x)-h(t-r,x), \ \forall \ x\in {\mathbb R}^d,
$$
and
$$
h^\delta(t-r,x_1,x_2)=h(t+\delta-r,x_1)-h(t+\delta-r,x_2), \ \forall \ x_1,x_2\in {\mathbb R}^d.
$$

We first calculate the term $I_1^\delta$:
\begin{eqnarray}\label{2.10}
|I_1^\delta(t)|&\leq& C \int_0^t\int_{|x-z|\leq 2|x-y|}[h^\delta(t-r)]_\alpha |x-z|^\alpha e^{-\frac{|x-z|^2}{2r}}r^{-\frac{d+1}{2}}dzdr
\nonumber\\ &\leq& C\|h^\delta\|_{L^q(0,t;\mathcal{C}^\alpha_b(\mathbb{R}^d))}\int_{|x-z|\leq 2|x-y|}|x-z|^\alpha\Big[ \int_0^te^{-\frac{q^\prime|x-z|^2}{2r}}
r^{-\frac{(d+1)q^\prime}{2}}dr\Big]^{\frac{1}{q^\prime}}dz
\nonumber\\ &\leq & C\|h^\delta\|_{L^q(0,t;\mathcal{C}^\alpha_b(\mathbb{R}^d))}\int_{|x-z|\leq 2|x-y|}|x-z|^{1+\alpha-\frac{2}{q}-d}dz\nonumber\\ &\leq & C\|h^\delta\|_{L^q(0,t;\mathcal{C}^\alpha_b(\mathbb{R}^d))}
|x-y|^{\theta},
\end{eqnarray}
where in the second line we have used the H\"{o}lder and Minkowski inequalities, and $1/q^\prime+1/q=1$.

Similarly, we get
\begin{eqnarray}\label{2.11}
|I_2^\delta(t)|\leq C\|h^\delta\|_{L^q(0,t;\mathcal{C}^\alpha_b(\mathbb{R}^d))}|x-y|^\theta, \ \ t\in [0,T), \ t+\delta\leq T.
\end{eqnarray}

For the term $I_3^\delta$, with the aid of Gauss-Green's formula, the H\"{o}lder and Minkowski inequalities, then
\begin{eqnarray}\label{2.12}
|I_3^\delta(t)|&=&\Big|\int_0^tdr\int_{|x-z|=2|x-y|}K(r,y-z)n_i
[h^\delta(t-r,y)-h^\delta(t-r,x)]dS\Big|
\nonumber\\&\leq & C\|h^\delta\|_{L^q(0,t;\mathcal{C}^\alpha_b(\mathbb{R}^d))}|x-y|^\alpha
\Big[\int_0^t\Big(\int_{|x-z|= 2|x-y|} r^{-\frac{d}{2}}e^{-\frac{|x-z|^2}{2r}} dS\Big)^{q^\prime}dr\Big]^{\frac{1}{q^\prime}}\nonumber\\
&\leq & C\|h^\delta\|_{L^q(0,t;\mathcal{C}^\alpha_b(\mathbb{R}^d))}|x-y|^\alpha\int_{|x-z|=2|x-y|} |y-z|^{-d-\frac{2}{q}+2}dz\nonumber\\
&\leq  & C\|h^\delta\|_{L^q(0,t;\mathcal{C}^\alpha_b(\mathbb{R}^d))}|x-y|^\theta.
\end{eqnarray}

To estimate $I_4^\delta$, using the H\"{o}lder inequality first, the Minkowski integral inequality next, then
\begin{eqnarray*}
|I_4^\delta(t)|\leq C\|h^\delta\|_{L^q(0,t;\mathcal{C}^\alpha_b(\mathbb{R}^d))} \int_{|x-z|> 2|x-y|}|x-z|^\alpha \Big(\int_0^t|\partial_{x_i}K(r,x-z)-
\partial_{y_i}K(r,y-z)|^{q^\prime}dr\Big)^{\frac{1}{q^\prime}}dz.
\end{eqnarray*}
For every $\eta\in [x,y]$, a segment of $x$ and $y$, due to $|x-z|>2|x-y|$, then
$$
\frac{1}{2}|x-z| \leq |\eta-z|\leq 2|x-z|.
$$
By virtue of mean value inequality and the property of second order partial
derivatives of the heat kernel $K(t,x)$, it yields that
\begin{eqnarray}\label{2.13}
|I_4^\delta(t)| &\leq & C\|h^\delta\|_{L^q(0,t;\mathcal{C}^\alpha_b(\mathbb{R}^d))} |x-y|\int_{|x-z|> 2|x-y|}|x-z|^\alpha\Big(\int_0^\infty r^{-\frac{(d+2)q^\prime}{2}}e^{-\frac{q^\prime|x-z|^2}{8r}}dr\Big)^{\frac{1}{q^\prime}}dz \nonumber\\
&\leq & C\|h^\delta\|_{L^q(0,t;\mathcal{C}^\alpha_b(\mathbb{R}^d))}|x-y|\int_{|x-z|> 2|x-y|}|x-z|^{\alpha-d-2+\frac{2}{q^\prime}} dz \nonumber\\&\leq& C\|h^\delta\|_{L^q(0,t;\mathcal{C}^\alpha_b(\mathbb{R}^d))}|x-y|^\theta.
\end{eqnarray}

Proceeding as the calculations from (\ref{2.10}) to (\ref{2.13}) lead to
\begin{eqnarray}\label{2.14}
\max\{|I_5^\delta(t)|,|I_6^\delta(t)|,|I_7^\delta(t)|,|I_8^\delta(t)|\}
\leq C\|h\|_{L^q(0,\delta;\mathcal{C}^\alpha_b(\mathbb{R}^d))}|x-y|^\theta, \ t,t+\delta\in [0,T].
\end{eqnarray}

Since $1\leq i\leq d$, combining (\ref{2.10}) to (\ref{2.14}), one asserts that
\begin{eqnarray}\label{2.15}
[\nabla u(t+\delta)-\nabla u(t)]_\theta\leq C [\|h^\delta\|_{L^q(0,t;\mathcal{C}^\alpha_b(\mathbb{R}^d))}+\|h\|_{L^q(0,\delta;\mathcal{C}^\alpha_b(\mathbb{R}^d))}], \ t,t+\delta\in [0,T].
\end{eqnarray}
Analogue calculations from (\ref{2.10}) to (\ref{2.14}) also hints
\begin{eqnarray}\label{2.16}
\|u(t+\delta)-u(t)\|_{\mathcal{C}^1_b(\mathbb{R}^d)} \leq C [\|h^\delta\|_{L^q(0,t;\mathcal{C}^\alpha_b(\mathbb{R}^d))}+\|h\|_{L^q(0,\delta;\mathcal{C}^\alpha_b(\mathbb{R}^d))}], \ t,t+\delta\in [0,T].
\end{eqnarray}
From (\ref{2.15}) and (\ref{2.16}), then
\begin{eqnarray}\label{2.17}
\|u(t+\delta)-u(t)\|_{\mathcal{C}^{1,\theta}_b(\mathbb{R}^d)}\leq C [\|h^\delta\|_{L^q(0,t;\mathcal{C}^\alpha_b(\mathbb{R}^d))}+\|h\|_{L^q(0,\delta;\mathcal{C}^\alpha_b(\mathbb{R}^d))}], \ t,t+\delta\in [0,T].
\end{eqnarray}
By letting $\delta$ tend to zero in (\ref{2.17}), we prove that as a ${\mathcal C}_b^{1,\theta}({\mathbb R}^d)$ valued function, $u$ is right continuous in $t$. If one replaces $\delta$ by $-\delta$, by repeating above calculations, we derive the left continuity of $u$ in $t$. In particular, if one chooses $t=0$ and $\delta=T$ in (\ref{2.17}) then (\ref{2.4}) holds.

$(ii)$ We set a mapping ${\mathcal T}$ on ${\mathcal C}([0,T];{\mathcal C}_b^{1,\alpha}({\mathbb R}^d))$ by
\begin{eqnarray}\label{2.18}
{\mathcal T} v(t,x)=\int_0^{t}K(r,\cdot)\ast (g(t-r,\cdot)\cdot \nabla v(t-r,\cdot))(x)dr+\int_0^tK(r,\cdot)\ast h(t-r,\cdot)(x)dr.
\end{eqnarray}
Then ${\mathcal T} v\in L^q(0,T;\mathcal{C}^{2,\alpha}_b(\mathbb{R}^d))
\cap W^{1,q}(0,T;\mathcal{C}^{\alpha}_b(\mathbb{R}^d))$. To finish the
proof, it suffices to show that the mapping is contractive on
${\mathcal C}([0,T];{\mathcal C}_b^{1,\alpha}({\mathbb R}^d))$, so there
is a unique $u\in{\mathcal C}([0,T];{\mathcal C}_b^{1,\alpha}({\mathbb R}^d))$
satisfying $u={\mathcal T} u$. This fact combining an argument as $g=0$
implies the existence of strong solutions of the Cauchy problem
(\ref{2.1}).

\smallskip
Let $t>0$ be given, $v_1,v_2\in {\mathcal C}([0,T];{\mathcal C}_b^{1,\alpha}({\mathbb R}^d))$, by utilising (\ref{2.4}), then
\begin{eqnarray*}
\|{\mathcal T} v_1-{\mathcal T} v_2\|_{{\mathcal C}([0,t];{\mathcal C}_b^{1,\alpha}({\mathbb R}^d))}
&\leq& C \||\nabla v_1-\nabla v_2|g\|_{L^2(0,t;\mathcal{C}^\alpha_b(\mathbb{R}^d))}\nonumber\\
&\leq & C\|g\|_{L^2(0,t;\mathcal{C}^\alpha_b(\mathbb{R}^d;\mR^d))}\|v_1-v_2\|_{{\mathcal C}([0,t];{\mathcal C}_b^{1,\alpha}({\mathbb R}^d))},
\end{eqnarray*}
which suggests that if $t>0$ is sufficiently small, then ${\mathcal T}$ is contractive, so there is a unique $u\in {\mathcal C}([0,t];{\mathcal C}_b^{1,\alpha}({\mathbb R}^d))$ such that ${\mathcal T} u=u$. For such small $t$, we assert that
\begin{eqnarray}\label{2.19}
\|u\|_{{\mathcal C}([0,t];{\mathcal C}_b^{1,\alpha}({\mathbb R}^d))}
&\leq& C [\|g\cdot \nabla u\|_{L^2(0,t;\mathcal{C}^\alpha_b(\mathbb{R}^d))}+\|h\|_{L^2(0,t;\mathcal{C}^\alpha_b(\mathbb{R}^d))}]\nonumber\\
&\leq & C\|g\|_{L^2(0,t;\mathcal{C}^\alpha_b(\mathbb{R}^d;\mR^d))}\|u\|_{{\mathcal C}([0,t];{\mathcal C}_b^{1,\alpha}({\mathbb R}^d))}+\|h\|_{L^2(0,t;\mathcal{C}^\alpha_b(\mathbb{R}^d))}].
\end{eqnarray}
Since now  $C\|g\|_{L^2(0,t;\mathcal{C}^\alpha_b(\mathbb{R}^d;\mR^d))}<1$, from (\ref{2.19}), then
\begin{eqnarray*}
\|u\|_{{\mathcal C}([0,t];{\mathcal C}_b^{1,\alpha}({\mathbb R}^d))}
\leq \frac{C}{1-C\|g\|_{L^2(0,t;\mathcal{C}^\alpha_b(\mathbb{R}^d;\mR^d))}}\|h\|_{L^2(0,t;\mathcal{C}^\alpha_b(\mathbb{R}^d))}\leq C\|h\|_{L^2(0,T;\mathcal{C}^\alpha_b(\mathbb{R}^d))}.
\end{eqnarray*}
We then repeat the proceeding arguments to extend the solution to the time interval $[t,2t]$. Continuing this procedure with finitely many steps, we construct a solution on $[0,T]$ for any given $T>0$ and get the inequality  (\ref{2.5}) on $[0,T]$.  The argument appeals above for $\|u\|_{{\mathcal C}([0,T];{\mathcal C}_b^{1,\alpha}({\mathbb R}^d))}$, adapted to $\|u\|_{L^q(0,T;\mathcal{C}^{2,\alpha}_b(\mathbb{R}^d))}$ now, yields the inequality (\ref{2.6}), and from this we complete the proof. $\Box$

From Lemma \ref{lem2.2}, we have
\begin{corollary} \label{cor2.1} Let $\alpha,h$ and $g$ be described in Lemma \ref{lem2.2} (ii). For $\lambda> 0$,
consider the Cauchy problem
\begin{eqnarray}\label{2.20}
\left\{\begin{array}{ll}
\partial_{t}u(t,x)=\frac{1}{2}\Delta u(t,x)+g(t,x)\cdot \nabla u(t,x)+h(t,x)-\lambda u(t,x), \ (t,x)\in (0,T)\times {\mathbb R}^d, \\
u(0,x)=0, \  x\in{\mathbb R}^d.
\end{array}\right.
\end{eqnarray}
Then there is a unique strong solution $u\in L^2(0,T;\mathcal{C}^{2,\alpha}_b(\mathbb{R}^d))\cap W^{1,2}(0,T;\mathcal{C}^{\alpha}_b(\mathbb{R}^d))$  to (\ref{2.20}). Moreover, there is a real number $\varepsilon>0$ such that
\begin{eqnarray}\label{2.21}
\|\nabla u\|_{{\mathcal C}([0,T];{\mathcal C}_b^0({\mathbb R}^d))}\leq C\lambda^{-\varepsilon}\|h\|_{L^2(0,T;\mathcal{C}^\alpha_b(\mathbb{R}^d))}
[\|h\|_{L^2(0,T;\mathcal{C}^\alpha_b(\mathbb{R}^d))}+\|g\|_{L^2(0,T;\mathcal{C}^\alpha_b(\mathbb{R}^d;\mR^d))}],
\end{eqnarray}
where the constant $C$ in (\ref{2.21}) only depends on $\alpha,d$ and $\|g\|_{L^2(0,T;\mathcal{C}^\alpha_b(\mathbb{R}^d;\mR^d))}$.
\end{corollary}
\vskip2mm\noindent
\textbf{Proof.} Let $v(t,x)=e^{\lambda t}u(t,x)$. If $u$ is a strong solution of (\ref{2.20}), then $v$ is a strong solution of the following Cauchy problem
\begin{eqnarray}\label{2.22}
\left\{\begin{array}{ll}
\partial_{t}v(t,x)=\frac{1}{2}\Delta v(t,x)+g(t,x)\cdot \nabla v(t,x)+e^{\lambda t} h(t,x), \ (t,x)\in (0,T)\times {\mathbb R}^d, \\
v(0,x)=0, \  x\in{\mathbb R}^d,
\end{array}\right.
\end{eqnarray}
and vice versa. By Lemma \ref{lem2.2}, (\ref{2.22}) exists a unique strong solution, and there is a constant $C>0$ such that
\begin{eqnarray*}
\|v\|_{{\mathcal C}([0,T];{\mathcal C}_b^{1,\alpha}({\mathbb R}^d))} \leq C\|e^{\lambda t}h\|_{L^2(0,T;\mathcal{C}^{\alpha}_b(\mathbb{R}^d))}.
\end{eqnarray*}
Therefore, (\ref{2.20}) exists a unique strong solution $u$, and there is a constant $C>0$ such that
\begin{eqnarray}\label{2.23}
\|u\|_{{\mathcal C}([0,T];{\mathcal C}_b^{1,\alpha}({\mathbb R}^d))}\leq C\| h\|_{L^2(0,T;\mathcal{C}^{\alpha}_b(\mathbb{R}^d))}.
\end{eqnarray}

By Lemma \ref{lem2.1}, the unique solution satisfies the following integral equation:
\begin{eqnarray*}
u(t,x)=\int_0^te^{-\lambda r} K(r,\cdot)\ast (g(t-r,\cdot)\cdot \nabla u(t-r,\cdot))(x)dr+\int_0^te^{-\lambda r} K(r,\cdot)\ast h(t-r,\cdot)dr.
\end{eqnarray*}
For every $x\in{\mathbb R}^d$ and $1\leq i\leq d$,
\begin{eqnarray*}
&&|\partial_{x_i}u(t,x)|\nonumber\\&\leq&\int_0^tdr\int_{\mathbb{R}^d}
e^{-\lambda r}|\partial_{x_i}K(r,x-z)||g(t-r,z)\cdot \nabla u(t-r,z)-g(t-r,x)\cdot \nabla u(t-r,x)|dz\nonumber\\&&+
\int_0^tdr\int_{\mathbb{R}^d}e^{-\lambda r}
|\partial_{x_i}K(r,x-z)||h(t-r,z)-h(t-r,x)|dz
\nonumber\\&\leq& \int_0^te^{-\lambda r}[g(t-r,\cdot)\cdot \nabla u(t-r,\cdot)]_\alpha dr\int_{\mathbb{R}^d} |x-z|^\alpha e^{-\frac{|x-z|^2}{2r}}r^{-\frac{d+1}{2}}dz
\nonumber\\&&+
\int_0^te^{-\lambda r}[h(t-r)]_\alpha dr\int_{\mathbb{R}^d} |x-z|^\alpha e^{-\frac{|x-z|^2}{2r}}r^{-\frac{d+1}{2}}dz
\nonumber\\&\leq&
C\|u\|_{{\mathcal C}([0,T];{\mathcal C}_b^{1,\alpha}({\mathbb R}^d))}\int_0^te^{-\lambda r}r^{-\frac{1-\alpha}{2}}
\|g(t-r)\|_{\mathcal{C}^\alpha_b(\mathbb{R}^d;\mR^d)} dr
+C\int_0^te^{-\lambda r}r^{-\frac{1-\alpha}{2}}
[h(t-r)]_\alpha dr
\nonumber\\&\leq & C\|h\|_{L^2(0,T;\mathcal{C}^\alpha_b(\mathbb{R}^d))}[\|h\|_{L^2(0,T;\mathcal{C}^\alpha_b(\mathbb{R}^d)}+
\|g\|_{L^2(0,T;\mathcal{C}^\alpha_b(\mathbb{R}^d;\mR^d))}]
\Big[\int_0^tr^{-\frac{(1-\alpha)p_1}{2}}dr\Big]^{\frac{1}{p_1}}
\Big[\int_0^te^{-\frac{\lambda r}{\varepsilon}}dr\Big]^{\varepsilon}
\nonumber\\&\leq & C
\lambda^{-\varepsilon},
\end{eqnarray*}
where $p_1=1/(1-\alpha)+1$, $\varepsilon=1/2-1/p_1$, and in the third inequality we have used (\ref{2.23}). $\Box$

We now extend the constant coefficients equation (\ref{2.20}) to a variable coefficients equation and found an analogue of Corollary \ref{cor2.1}. To be precise, we consider  the following Cauchy problem
\begin{eqnarray}\label{2.24}
\left\{\begin{array}{ll}
\partial_{t}u(t,x)=\frac{1}{2}\sum_{i,j=1}^da_{i,j}(t,x)\partial^2_{x_i,x_j} u(t,x)+g(t,x)\cdot\nabla u(t,x)\\ \qquad\qquad \ \ +h(t,x)-\lambda u(t,x), \ \ (t,x)\in (0,T)\times {\mathbb R}^d, \\
u(0,x)=0, \  x\in{\mathbb R}^d,
\end{array}\right.
\end{eqnarray}
where $a_{i,j}(t,x), i,j=1,\ldots,d$ are real-valued functions such that $a_{i,j}\in L^\infty(0,T;\mathcal{C}^\alpha_b(\mathbb{R}^d))$ and for every $(t,x)\in (0,T)\times \mathbb{R}^d$, (\ref{1.6}) holds.

From \cite[Theorem 2.3]{Chen17} or \cite[Theorem 1.8]{Chen18}, there is a unique continuous function $p(t,x)$ called the fundamental solution or heat kernel of the differential operator $\frac{1}{2}\sum_{i,j=1}^da_{i,j}(t,x)\partial^2_{x_i,x_j}$, with the following properties:

$(H_1)$ (Two-sides estimates) For every $t\in (0,T]$ and $x\in \mR^d$, there exist constants $C,\mu>0$ such that
\begin{eqnarray}\label{2.25}
C^{-1}t^{-\frac{d}{2}}e^{-\frac{|x|^2}{t\mu}}\leq p(t,x)\leq Ct^{-\frac{d}{2}}e^{-\frac{\mu|x|^2}{t}}
\end{eqnarray}

$(H_2)$ (Gradient estimate) For every $t\in (0,T]$ and $x\in \mR^d$, there exist constants $C,\mu>0$ such that for $m=0,1,2$
\begin{eqnarray}\label{2.26}
|\nabla^mp(t,x)|\leq Ct^{-\frac{d+m}{2}}e^{-\frac{\mu |x|^2}{t}}.
\end{eqnarray}
By (\ref{2.25}) and (\ref{2.26}), then \cite[Theorem 3.3]{Kry02} is applicable for the Cauchy problem
\begin{eqnarray}\label{2.27}
\left\{\begin{array}{ll}
\partial_{t}u(t,x)=\frac{1}{2}\sum_{i,j=1}^da_{i,j}(t,x)\partial^2_{x_i,x_j} u(t,x) +h(t,x), \ \ (t,x)\in (0,T)\times {\mathbb R}^d, \\
u(0,x)=0, \  x\in{\mathbb R}^d.
\end{array}\right.
\end{eqnarray}
For $h\in L^2(0,T;\mathcal{C}^{\alpha}_b(\mathbb{R}^d))$, there is a unique $u\in L^2(0,T;\mathcal{C}^{2,\alpha}_b(\mathbb{R}^d))\cap W^{1,2}(0,T;\mathcal{C}^{\alpha}_b(\mathbb{R}^d))$ solving the Cauchy problem
(\ref{2.27}). The same arguments for Lemmas \ref{lem2.1}, \ref{lem2.2}, and Corollary \ref{cor2.1} used again, we conclude that

\begin{theorem} \label{thm2.1} Let $\alpha\in (0,1)$. Suppose that $g\in L^q(0,T;\mathcal{C}^{\alpha}_b(\mathbb{R}^d;\mathbb{R}^d))$ and $h\in L^q(0,T;\mathcal{C}^{\alpha}_b(\mathbb{R}^d))$. Let $(a_{i,j})$ be a
symmetric $d\times d$ matrix valued function whose components
$a_{i,j}$ are in $L^\infty(0,T;\mathcal{C}^{\alpha}_b(\mathbb{R}^d))$, and let
(\ref{1.6}) hold. Then
there is a unique strong solution $u\in L^2(0,T;\mathcal{C}^{2,\alpha}_b(\mathbb{R}^d))\cap W^{1,2}(0,T;\mathcal{C}^{\alpha}_b(\mathbb{R}^d))$ to (\ref{2.24}).
Moreover, $u\in {\mathcal C}([0,T];{\mathcal C}_b^{1,\alpha}({\mathbb R}^d))$ and for all $\lambda>0$, (\ref{2.21}) holds.
\end{theorem}

\begin{remark}\label{rem2.1} For the Cauchy problem (\ref{2.24}), when
$g\in L^\infty(0,T;\mathcal{C}^{\alpha}_b(\mathbb{R}^d;\mathbb{R}^d))$,
the existence and unique of strong solution has been proved by Krylov
\cite{Kry02}, and when
$g\in B([0,T];\mathcal{C}^{\alpha}_b(\mathbb{R}^d;\mathbb{R}^d))$,
$$
\sup_{t\in[0,T]}\|g(t)\|_{\mathcal{C}^{\alpha}_b(\mathbb{R}^d;\mathbb{R}^d)}<\infty,
$$
the existence and unique of
$B([0,T];\mathcal{C}^{2,\alpha}_b(\mathbb{R}^d))$ solution
is also established by Lorenzi \cite{Lor}.
Here we only assume that
$g\in L^2(0,T;\mathcal{C}_b^\alpha(\mathbb{R}^d;\mathbb{R}^d))$,
so we extend Krylov and Lorenzi's results. This result plays a central role
in proving the pathwise uniqueness of solutions since it yields the following It\^{o} formula.
\end{remark}

\begin{theorem} \label{thm2.2} Let $X_t$ satisfy (\ref{1.1}) with
$b\in L^2(0,T;\mathcal{C}^{\alpha}_b(\mathbb{R}^d;\mathbb{R}^d))$ and
$\sigma_{i,j}\in L^\infty(0,T;\mathcal{C}^{\alpha}_b(\mathbb{R}^d))$.
Let $\alpha$ and $(a_{i,j})=
(\sigma_{i,k}\sigma_{j,k})$ be described in Theorem \ref{thm2.1}, and let $u\in L^2(0,T;\mathcal{C}^{2,\alpha}_b(\mathbb{R}^d))\cap
W^{1,2}(0,T;\mathcal{C}^{\alpha}_b(\mathbb{R}^d))$. Then the
following It\^{o} formula
\begin{eqnarray}\label{2.28}
u(t,X_t(x))&=&u(0,x)+\int_0^t[\partial_ru(r,X_r(x))+
\nabla u(r,X_r(x))\cdot b(r,X_r(x))]dr
\nonumber\\&&+\frac{1}{2}\sum_{i,j=1}^d\int_0^ta_{i,j}(r,X_r(x))\partial^2_{X_i,X_j}u(r,X_r(x))dr
\nonumber\\&&+
\sum_{i,k=1}^d\int_0^t\partial_{X_i}u(r,X_r(x))\sigma_{i,k}(r,X_r(x))dW_{k,r}
\end{eqnarray}
holds, for every $t\in [0,T]$.
\end{theorem}
\vskip0mm\noindent
\textbf{Proof.} For $0<\varepsilon<1$,
a given function $h\in L^2(0,T;\mathcal{C}^\alpha_b(\mathbb{R}^d))$,
and every $t\in [0,T]$, we set
\begin{eqnarray*}
h_\varepsilon(t,x)=\int_0^1h(t+r\varepsilon,x)dr.
\end{eqnarray*}
By the Lebesgue differentiation theorem, the random
variable $h_\varepsilon(t,X_t(x))$ converges
to $h(t,X_t(x))$ ${\mathbb P}-a.s.$ as $\varepsilon\downarrow 0$
for almost
every $t\in [0,T]$. Let $u$ be stated in Theorem \ref{thm2.2} and we define $u_\varepsilon(t,x)$ as the previous procedure.  Then
$u_\varepsilon(t,x)$ is continuous and differentiable in $t$ and
the random variable $u_\varepsilon(t,X_t(x))$ converges to $u(t,X_t(x))$ almost
surely as $\varepsilon\downarrow 0$. We apply the classical
It\^{o} formula to $u_\varepsilon(t,X_t)$ and get
\begin{eqnarray}\label{2.29}
u_\varepsilon(t,X_t(x))&=&u_\varepsilon(0,x)+\int_0^t[\partial_su_\varepsilon(r,X_r(x))+
\nabla u_\varepsilon(r,X_r(x))\cdot b(r,X_r(x))]dr
\nonumber\\&&+\frac{1}{2}\sum_{i,j=1}^d\int_0^ta_{i,j}(r,X_r(x))\partial^2_{x_i,x_j}u_\varepsilon(r,X_r(x))dr
\nonumber\\&&+
\sum_{i,k=1}^d\int_0^t\partial_{X_i}u_\varepsilon(r,X_r(x))\sigma_{i,k}(r,X_r(x))dW_{k,r}.
\end{eqnarray}
For every given function $h\in L^2(0,T;\mathcal{C}^\alpha_b(\mathbb{R}^d))$,
we estimate $h_\varepsilon(t,X_t(x))$ by
\begin{eqnarray*}
|h_\varepsilon(t,X_t(x))|\leq \int_0^1\|h(t+r\varepsilon)\|_{{\mathcal C}_b({\mathbb R}^d)}dr
\leq
\sup_{0<\varepsilon<1}\frac{1}{\varepsilon}
\int_t^{t+\varepsilon}\|h(r)\|_{{\mathcal C}_b({\mathbb R}^d)}dr=:g(t).
\end{eqnarray*}
With the aid the property for Hardy-Littlewood maximum function \cite[Theorem 1, pp.5]{Stein}, then $g\in L^2(0,T)$.
By applying the dominated convergence theorem, random variable
$\int_0^t\zeta(r)h_\varepsilon(r,X_r(x))dr$
converges to $\int_0^t\zeta(r)h(r,X_r(x))dr$ ${\mathbb P}-a.s.$ as
$\varepsilon \downarrow0$ for every $\zeta\in L^2(0,T)$, which suggests that
the second and third terms in the right hand side of (\ref{2.29}) converge to
the second and third terms in the right hand side of (\ref{2.28}), respectively.

We calculate the difference for the last terms in (\ref{2.29}) and (\ref{2.28}) by
\begin{eqnarray}\label{2.30}
&&{\mathbb E}\Big|\sum_{i,k=1}^d\int_0^t[\partial_{X_i}u_\varepsilon(r,X_r(x))-\partial_{X_i}u(r,X_r(x))]\sigma_{i,k}(r,X_r(x))
dW_{k,r}\Big|\nonumber\\&=&\sum_{i,k=1}^d{\mathbb E}\int_0^t|[\partial_{X_i}u_\varepsilon(r,X_r(x))-
\partial_{X_i}u(r,X_r(x))]\sigma_{i,k}(r,X_r(x))|^2dr
\nonumber\\&\leq &C{\mathbb E}\int_0^t|\nabla u_\varepsilon(r,X_r(x))-
\nabla u(r,X_r(x))|^2dr.
\end{eqnarray}
Clearly, $\nabla u_\varepsilon(r,X_r(x))-
\nabla u(r,X_r(x))$ converges
to 0, ${\mathbb P}-a.s.$ as $\varepsilon\downarrow 0$ for almost
every $r\in [0,T]$, and
\begin{eqnarray*}
|\nabla u_\varepsilon(r,X_r(x))-
\nabla u(r,X_r(x))|
\leq \sup_{0<\varepsilon<1}\frac{1}{\varepsilon}
\int_r^{r+\varepsilon}\|\nabla u(s)\|_{{\mathcal C}_b({\mathbb R}^d)}ds
+\|\nabla u(r)\|_{{\mathcal C}_b({\mathbb R}^d)}\in L^2(0,T).
\end{eqnarray*}
Thus, the last term in the right hand side of (\ref{2.30}) vanishes
if $\varepsilon$
tends to $0$, and it also suggests that the random variable
$\int_0^t\partial_{X_i}u_\varepsilon(r,X_r(x))\sigma_{i,k}(r,X_r(x))dW_{k,r}$
converges to the random variable $\int_0^t\partial_{X_i}u(r,X_r(x))\sigma_{i,k}(r,X_r(x))dW_{k,r}$
${\mathbb P}-a.s.$ up to choosing a unlabelled subsequence. $\Box$

\section{Proof of Theorem \ref{thm1.1}}
\label{sec3}\setcounter{equation}{0}
Firstly, we present an approximating result.
\begin{lemma} \label{lem3.1} (i) Let $T>0$ be a real number, and let
$h\in L^q(0,T;\mathcal{C}_{u}(\mathbb{R}^d))$ with $q\in[1,2]$. We set $h^n(t,x)=(h(t,\cdot)\ast\rho_n)(x), n\in\mathbb{N}$, where
$\ast$ stands for the usual convolution and $\rho_n(x)=n^d \rho(nx)$ with
\begin{eqnarray}\label{3.1}
0\leq \rho \in \mathcal{C}^\infty_0(\mathbb{R}^d) , \ \ \mbox{support}(\rho)\subset B_0(1), \ \int_{\mathbb{R}^d}\rho(x)dx=1.
\end{eqnarray}
Then
\begin{eqnarray}\label{3.2}
\lim_{n\rightarrow \infty}\|h^n-h\|_{L^1(0,T;\mathcal{C}_{u}(\mathbb{R}^d))}
=\lim_{n\rightarrow \infty}\int_0^T\sup_{x\in\mathbb{R}^d}|h^n(t,x)-h(t,x)|^qdt=0.
\end{eqnarray}

(ii) Let $W_t,W_t^n,n=1,2,\ldots$ be $d$-dimensional standard Wiener
processes on a same stochastic basis
$(\Omega,{\mathcal F},\{{\mathcal F}_t\}_{0\leq t\leq T},{\mathbb P})$
for which $W_\cdot^n$ converges to $W_\cdot$, ${\mathbb P}$-a.s..
Assume $f\in L^2(0,T)$, then
\begin{eqnarray}\label{3.3}
\lim_{n\rightarrow \infty}\mathbb{E}\Big|\int_0^Tf(t)d[W_t^n-W_t]\Big|^2=0.
\end{eqnarray}
\end{lemma}
\smallskip
\textbf{Proof.} $(i)$ For every $n\in {\mathbb N}$ and $t\in [0,T]$,
\begin{eqnarray*}
h^n(t,x)-h(t,x)=\int_{\mathbb{R}^d}[h(t,x-\frac{y}{n})-h(t,x)]\rho(y)dy,
\end{eqnarray*}
which suggests that
\begin{eqnarray*}
\sup_{x\in\mathbb{R}^d}|h^n(t,x)-h(t,x)|^q\leq
C\int_{|y|\leq 1}\sup_{x\in\mathbb{R}^d}|h(t,x-\frac{y}{n})-h(t,x)|^q\rho(y)dy.
\end{eqnarray*}
Therefore, (\ref{3.2}) is true.

\smallskip
$(ii)$ Since $f\in L^2(0,T)$, we approximate it by a sequence of smooth functions
$f_\varepsilon\in W^{1,2}(0,T)$ such that $f_\varepsilon(T)=0$ and
$f_\varepsilon\rightarrow f$ in $L^2(0,T)$ as $\varepsilon\downarrow0$. Then
\begin{eqnarray}\label{3.4}
&&\lim_{n\rightarrow \infty}\mathbb{E}\Big|\int_0^Tf(t)d[W_t^n-W_t]\Big|^2\nonumber\\
&\leq& 2\lim_{n\rightarrow \infty}\mathbb{E}\Big|\int_0^T[f(t)-f_\varepsilon(t)]d[W_t^n-W_t]\Big|^2
+2\lim_{n\rightarrow \infty}\mathbb{E}\Big|\int_0^Tf_\varepsilon(t)d[W_t^n-W_t]\Big|^2.
\end{eqnarray}

In view of It\^{o}'s isometry and the integration by parts for Wiener's
integral, from (\ref{3.4}), it yields that
\begin{eqnarray*}
&&\lim_{n\rightarrow \infty}\mathbb{E}\Big|\int_0^Tf(t)d[W_t^n-W_t]\Big|^2\nonumber\\&\leq& C\lim_{n\rightarrow \infty}\mathbb{E}\int_0^T|f(t)-f_\varepsilon(t)|^2dt
+2\lim_{n\rightarrow \infty}\mathbb{E}\Big|\int_0^Tf^\prime_\varepsilon(t)[W_t^n-W_t]dt\Big|^2
\nonumber\\&=& C\mathbb{E}\int_0^T|f(t)-f_\varepsilon(t)|^2dt.
\end{eqnarray*}
So (\ref{3.3}) holds by letting $\varepsilon\downarrow0$. $\Box$

In preparation to prove the H\"{o}lder continuity of solutions, we need another lemma.
\begin{lemma}(\cite[Exercise 2.10, pp.31]{RY} or \cite[Theorem 1.1]{WL}) \label{lem3.2} Let $\{X_t(x), x\in [0,1]^d, t\in [0,1]\}$ be a stochastic field for which there exist three strictly positive constants $\gamma,C, \varepsilon$ such that
\begin{eqnarray}\label{3.5}
\mathbb{E}[\sup_{0\leq t\leq 1}|X_t(x)-X_t(y)|^\gamma]\leq C|x-y|^{d+\varepsilon}.
\end{eqnarray}
Then there is a modification $\tilde{X}$ of $X$ such that
\begin{eqnarray}\label{3.6}
\mathbb{E}\Big[\sup_{0\leq t\leq 1}\Big(\sup_{x\neq y}\frac{|\tilde{X}_t(x)-\tilde{X}_t(y)|}{|x-y|^\beta}\Big)^\gamma\Big]<\infty
\end{eqnarray}
for every $\beta\in [0,\varepsilon/\gamma)$. In particular, the paths of $\tilde{X}$ are H\"{o}lder continuous in $x$ of order $\beta$.
\end{lemma}

We are now in a position to prove our main result to
SDE (\ref{1.1}).

$(i)$ We smooth out $b$ and $\sigma$ using the convolution:
$b^n(t,x)=(b(t,\cdot)\ast \rho_n)(x)$,
$\sigma^n(t,x)=(\sigma(t,\cdot)\ast \rho_n)(x)$ with $\rho_n$ given by
(\ref{3.1}). According to (\ref{3.2}), as $n\rightarrow \infty$,
\begin{eqnarray}\label{3.7}
\|b^n-b\|_{L^1(0,T;\mathcal{C}_{u}(\mathbb{R}^d))} \rightarrow 0, \ \ \|\sigma^n-\sigma\|_{L^2(0,T;\mathcal{C}_{u}(\mathbb{R}^d))} \rightarrow 0.
\end{eqnarray}
Moreover, for every $n\geq 1$, and almost every $t\in [0,T]$,
\begin{eqnarray}\label{3.8}
\|b^n(t)\|_{\mathcal{C}_{u}(\mathbb{R}^d)}\leq  \|b(t)\|_{\mathcal{C}_{u}(\mathbb{R}^d)}, \ \|\sigma^n(t)\|_{\mathcal{C}_{u}(\mathbb{R}^d)}\leq  \|\sigma(t)\|_{\mathcal{C}_{u}(\mathbb{R}^d)} .
\end{eqnarray}
Furthermore, there are two sequences $\theta^n$
and $l^n$, which are integrable and square-integrable functions on $[0,T]$, respectively,  such that
\begin{eqnarray*}
|b^n(t,x)-b^n(t,y)|\leq \theta^n(t)|x-y|, \ \ \forall \ x,y\in\mathbb{R}^d
\end{eqnarray*}
and
\begin{eqnarray*}
|\sigma^n(t,x)-\sigma^n(t,y)|\leq l^n(t)|x-y|, \ \ \forall \ x,y\in\mathbb{R}^d.
\end{eqnarray*}

By Cauchy-Lipschitz's theorem, there is a unique
$\{\mathcal{F}_t\}$-adapted, continuous, $d$-dimensional process $X_t^n$
defined for $t\in [0,T]$ on $(\Omega,\mathcal{F},\{\mathcal{F}_t\}_{0\leq t\leq T},\mathbb{P})$ such that
\begin{eqnarray}\label{3.9}
X_t^n=x+\int_0^tb^n(r,X^n_r)dr+\int_0^t\sigma^n(r,X^n_r)dW_r.
\end{eqnarray}

With the help of (\ref{3.8}), for every $0\leq t_1<t_2\leq T$,
\begin{eqnarray}\label{3.10}
\sup_{n}\mathbb{E}\int_{t_1}^{t_2}|b^n(t,X^n_t)|dt\leq  \int_{t_1}^{t_2}\sup_{x\in\mathbb{R}^d}|b(t,x)|dt
\end{eqnarray}
and
\begin{eqnarray}\label{3.11}
\sup_{n}\mathbb{E}\Big|\int_{t_1}^{t_2}\sigma^n(t,X^n_t)dW_t\Big|^2\leq \int_{t_1}^{t_2}\sup_{x\in\mathbb{R}^d}|\sigma(t,x)|^2dt.
\end{eqnarray}

Combining (\ref{3.9}), (\ref{3.10}) and (\ref{3.11}), for every $\epsilon>0$, one concludes that
\begin{eqnarray}\label{3.12}
\lim_{c\rightarrow \infty}\sup_n\sup_{0\leq t\leq T}\mathbb{P}\{|X_t^n|>c\}=0
\end{eqnarray}
and
\begin{eqnarray}\label{3.13}
\lim_{h\downarrow 0}\sup_{n\geq 1}
\sup_{0\leq t_1,t_2\leq T\atop{|t_1-t_2|\leq h}}\mathbb{P}\{|X_{t_1}^n-X_{t_2}^n|>\epsilon\}=0.
\end{eqnarray}

From (\ref{3.12}) and (\ref{3.13}), along with  Prohorov's theorem,
there is a subsequence still denoted by itself such that
$(X_\cdot^n,W_\cdot)$ weakly converge. Next, Skorohod's
representation theorem implies that there is a probability space
$(\tilde{\Omega},\tilde{\mathcal{F}},\{\tilde{\mathcal{F}}_t\}_{0\leq t\leq T},\tilde{\mathbb{P}})$
and random processes $(\tilde{X}_t^n,\tilde{W}_t^n)$, $(\tilde{X}_t,\tilde{W}_t)$
on this probability space such that

$(1)$ the finite dimensional distributions of $(\tilde{X}_t^n,\tilde{W}_t^n)$
coincide with the corresponding finite dimensional distributions of $(X_t^n,W_t)$.

$(2)$ $(\tilde{X}_\cdot^n,\tilde{W}_\cdot^n)$  converges to $(\tilde{X}_\cdot,\tilde{W}_\cdot)$, $\tilde{\mathbb{P}}-a.s.$.

In particular, $\tilde{W}$ is still a Wiener process and
\begin{eqnarray}\label{3.14}
\tilde{X}_t^n=x+\int_0^tb^n(r,\tilde{X}^n_r)dr+\int_0^t\sigma^n
(r,\tilde{X}^n_r)d\tilde{W}^n_r.
\end{eqnarray}
For $k\in{\mathbb N}$, then
\begin{eqnarray}\label{3.15}
&&\tilde{\mathbb{E}}\Big(\int_0^T|b^n(r,\tilde{X}^n_r)-
b(r,\tilde{X}_r)|dr\Big)\nonumber\\&\leq& \tilde{\mathbb{E}}\Big(\int_0^T|b^n(r,\tilde{X}^n_r)-
b^k(r,\tilde{X}^n_r)|dr\Big)
+\tilde{\mathbb{E}}\Big(\int_0^T|b^k(r,\tilde{X}^n_r)-
b^k(r,\tilde{X}_r)|dr\Big)
\nonumber\\&&+\tilde{\mathbb{E}}\Big(\int_0^T|b^k(r,\tilde{X}_r)-
b(r,\tilde{X}_r)|dr\Big)
\nonumber\\&\leq& C\Big[ \|b^n-b^k\|_{L^1(0,T;\mathcal{C}_{u}(\mathbb{R}^d))}+
\|b^k-b\|_{L^1(0,T;\mathcal{C}_{u}(\mathbb{R}^d))}\Big]
\nonumber\\&&+
\tilde{\mathbb{E}}\Big(\int_0^T|b^k(r,\tilde{X}^n_r)-b^k(r,\tilde{X}_r)|dr\Big).
\end{eqnarray}
We approach $n\rightarrow \infty$ first, $k\rightarrow \infty$ next,
from (\ref{3.7}) and (\ref{3.15}), it follows that
\begin{eqnarray}\label{3.16}
\lim_{n\rightarrow \infty}\int_0^tb^n(r,\tilde{X}^n_r)dr=\int_0^tb(r,\tilde{X}_r)dr,\,\, \tilde{\mathbb{P}}-a.s..
\end{eqnarray}

Similar calculations also suggest that
\begin{eqnarray*}
&&\tilde{\mathbb{E}}\Big|\int_0^T\sigma^n(r,\tilde{X}^n_r)d
\tilde{W}^n_r-
\int_0^T\sigma(r,\tilde{X}_r)d\tilde{W}_r\Big|^2\nonumber\\&\leq&
C\tilde{\mathbb{E}}\Big|\int_0^T\sigma^n(r,\tilde{X}^n_r)
d\tilde{W}^n_r-
\int_0^T\sigma(r,\tilde{X}_r)d\tilde{W}_r^n\Big|^2
+C\tilde{\mathbb{E}}\Big|\int_0^T\sigma(r,\tilde{X}_r)d\tilde{W}_r^n-
\int_0^T\sigma(r,\tilde{X}_r)d\tilde{W}_r\Big|^2
\nonumber\\&\leq& C\tilde{\mathbb{E}}\int_0^T\Big|\sigma^n(r,\tilde{X}^n_r)-
\sigma(r,\tilde{X}_r)\Big|^2dr+C\tilde{\mathbb{E}}\Big|\int_0^T\sigma(r,\tilde{X}_r)d\tilde{W}_r^n-
\int_0^T\sigma(r,\tilde{X}_r)d\tilde{W}_r\Big|^2
\nonumber\\&=&:J_1^n+J_2^n.
\end{eqnarray*}
We adopt the same procedure as in (\ref{3.15}) to assert that $J_1^n\rightarrow 0$ as $n\rightarrow \infty$.

On the other hand, thanks to the definition of stochastic integral, then
\begin{eqnarray}\label{3.17}
&&\tilde{\mathbb{E}}\Big|\int_0^t\sigma(r,\tilde{X}_r)d[\tilde{W}^n_r-\tilde{W}_r]
\Big|^2\nonumber\\&=&
\lim_{k\rightarrow \infty}\tilde{\mathbb{E}}\Big|\sum_{i=1}^k\sigma(r_i\wedge t,\tilde{X}_{r_i})[\tilde{W}^n_{r_{i+1}}-\tilde{W}_{r_{i+1}}-
\tilde{W}^n_{r_i}+\tilde{W}_{r_i}]\Big|^2
\nonumber\\&=&
\lim_{k\rightarrow \infty}\tilde{\mathbb{E}}\sum_{i=1}^k|\sigma(r_i\wedge t,\tilde{X}_{r_i})|^2|\tilde{W}^n_{r_{i+1}}-\tilde{W}_{r_{i+1}}-
\tilde{W}^n_{r_i}+\tilde{W}_{r_i}|^2
\nonumber\\&\leq&\lim_{k\rightarrow \infty}\tilde{\mathbb{E}}\sum_{i=1}^k\sup_{x\in\mathbb{R}^d}|\sigma(r_i\wedge t,x)|^2|\tilde{W}^n_{r_{i+1}}-\tilde{W}_{r_{i+1}}-
\tilde{W}^n_{r_i}+\tilde{W}_{r_i}|^2
\nonumber\\&=&\tilde{\mathbb{E}}\Big|\int_0^t\sup_{x\in\mathbb{R}^d}|\sigma(r,x)|d[\tilde{W}^n_r-\tilde{W}_r]
\Big|^2.
\end{eqnarray}
According to (\ref{3.3}), from  (\ref{3.17}), so $J_2^n\rightarrow 0$ as
$n\rightarrow \infty$.  Therefore,
\begin{eqnarray}\label{3.18}
\lim_{n\rightarrow \infty}\int_0^t\sigma^n(r,\tilde{X}^n_r)d\tilde{W}^n_r=
\int_0^t\sigma(r,\tilde{X}_r)d\tilde{W}_r, \,\, \tilde{\mathbb{P}}-a.s..
\end{eqnarray}
Combining (\ref{3.14}), (\ref{3.16}) and (\ref{3.18}), then (\ref{1.5}) holds.

$(ii)$ We proceed to show the pathwise uniqueness by using It\^{o}-Tanaka's trick. Consider the following vector valued Cauchy problem on $(0,T)\times \mathbb{R}^d$,
\begin{eqnarray}\label{3.19}
&&\partial_tU(t,x)=\sum_{i,j=1}^d\frac{1}{2}a_{i,j}(T-t,x)\partial^2_{x_i,x_j} U(t,x)+\sum_{i=1}^db_i(T-t,x)\partial_{x_i} U(t,x)\nonumber\\&&\qquad\qquad\quad-\lambda U(t,x)+ b(T-t,x),
\end{eqnarray}
with initial data $U(0,x)=0$. In view of Theorem \ref{thm2.1}, there is a unique
$U\in L^2(0,T;\mathcal{C}^{2,\alpha}_b(\mathbb{R}^d))\cap W^{1,2}(0,T;\mathcal{C}^{\alpha}_b(\mathbb{R}^d))$
solving (\ref{3.19}). Moreover,
$U\in {\mathcal C}([0,T];{\mathcal C}_b^1({\mathbb R}^d;\mathbb{R}^d))$ and there is a real number $\varepsilon>0$ such that
\begin{eqnarray*}
\|U\|_{{\mathcal C}([0,T];{\mathcal C}_b^1({\mathbb R}^d))}\leq
C\lambda^{-\varepsilon}.
\end{eqnarray*}
So,
\begin{eqnarray*}
\|U\|_{{\mathcal C}([0,T];{\mathcal C}_b^1({\mathbb R}^d))}<\frac{1}{2},
\quad \mbox{if} \quad
\lambda>(2C)^{\frac{1}{\varepsilon}}.
\end{eqnarray*}
Let $\lambda$ be big enough ($\lambda>(2C)^{1/\varepsilon}$) and fixed. We define
\begin{eqnarray}\label{3.20}
\Phi(t,x)=x+U(T-t,x).
\end{eqnarray}
Obviously, $\Phi$ forms a non-singular diffeomorphism of class ${\mathcal C}^1$
uniformly in $t\in [0,T]$ and
\begin{eqnarray}\label{3.21}
\frac{1}{2}<\|\nabla\Phi\|_{\mathcal{C}([0,T];\mathcal{C}_b({\mathbb R}^d))}<\frac{3}{2},
\quad  \frac{2}{3}<\|\nabla\Psi\|_{\mathcal{C}([0,T];\mathcal{C}_b({\mathbb R}^d))}<2,
\end{eqnarray}
where $\Psi(t,\cdot)=\Phi^{-1}(t,\cdot)$. Moreover, as a measurable measurable function $\Phi(t,x)-x$ in $(t,x)$, it belongs to $L^2(0,T;\mathcal{C}^{2,\alpha}_b(\mathbb{R}^d;\mathbb{R}^d))
\cap W^{1,2}(0,T;\mathcal{C}^{\alpha}_b(\mathbb{R}^d;\mathbb{R}^d))$.

\smallskip
As a result of It\^{o}'s formula (Theorem \ref{thm2.2}), we assert
\begin{eqnarray*}
d\Phi(t,X_t)&=&-\partial_t u(T-t,X_t)dt+\sum_{i=1}^db_i(t,X_t)\partial_{x_i} U(T-t,X_t)dt
\nonumber\\&&+\frac{1}{2}\sum_{i,j=1}^da_{i,j}(t,x)\partial^2_{x_i,x_j} u(T-t,X_t)dt
\nonumber\\&&+\sum_{i,k=1}^d\partial_{x_i} U(T-t,X_t)\sigma_{i,k}dW_{k,t}+b(t,X_t)dt+\sigma(t,X_t)dW_t
\nonumber\\&=&(\nabla U(T-t,X_t)+I)\sigma(t,X_t)dW_t+\lambda U(T-t,X_t)dt.
\end{eqnarray*}
Denote $Y_t=\Phi(t,X_t)$, it yields that
\begin{eqnarray}\label{3.22}
dY_t&=&\lambda U(T-t,\Psi(t,Y_t))dt+(I+\nabla U(T-t,\Psi(t,Y_t))\sigma(t,\Psi(t,Y_t))dW_t\nonumber\\&=&:
\tilde{b}(t,Y_t)dt+\tilde{\sigma}(t,Y_t)dW_t,
\end{eqnarray}
with $Y_0=y=\Phi(0,x)$. Then $\tilde{b}\in \mathcal{C}([0,T];\mathcal{C}_b^1(\mathbb{R}^d;\mathbb{R}^d))$. On the other hand, for $1\leq i,j\leq d$, $\sigma_{i,j}\in L^\infty(0,T;\mathcal{C}_b^\alpha(\mathbb{R}^d))$ and $\nabla\sigma_{i,j}\in L^2(0,T;L^\infty(\mathbb{R}^d;\mR^d))$, thus $\tilde{\sigma}\in L^2(0,T;W^{1,\infty}(\mathbb{R}^d))$. Let $Y^1_t$ and $Y_t^2$ be two solutions of (\ref{3.22}). Using the It\^{o} formula to $|Y_t^1-Y_t^2|^4$, then
\begin{eqnarray}\label{3.23}
d|Y_t^1-Y_t^2|^4&\leq & 4\langle |Y_t^1-Y_t^2|^2(Y_t^1-Y_t^2), \tilde{b}(t,Y_t^1)-\tilde{b}(t,Y_t^2) \rangle  dt \nonumber\\&&+6|Y_t^1-Y_t^2|^2\|\tilde{\sigma}(t,Y_t^1)-\tilde{\sigma}(t,Y_t^2)\|^2dt \nonumber\\&& +4\langle |Y_t^1-Y_t^2|^2(Y_t^1-Y_t^2), (\tilde{\sigma}(t,Y_t^1)-\tilde{\sigma}(t,Y_t^2))dW_t \rangle
\nonumber\\&\leq & C(1+\kappa(t))|Y_t^1-Y_t^2|^4dt \nonumber\\&&+4\langle |Y_t^1-Y_t^2|^2(Y_t^1-Y_t^2), (\tilde{\sigma}(t,Y_t^1)-\tilde{\sigma}(t,Y_t^2))dW_t \rangle,
\end{eqnarray}
where $\kappa\in L^1(0,T)$.

By taking the expectation in (\ref{3.23}), it follows that
\begin{eqnarray*}
\mE|Y_t^1-Y_t^2|^4\leq  C\int_0^t(1+\kappa(r))\mE|Y_r^1-Y_r^2|^4dr,
\end{eqnarray*}
which suggests that
\begin{eqnarray}\label{3.24}
\sup_{0\leq t \leq T}\mE|Y_t^1-Y_t^2|^4=0,
\end{eqnarray}
if one uses the Gr\"{o}nwall inequality.

We now use the It\^{o} formula to $|Y_t^1-Y_t^2|^2$ again first, BDG inequality next, H\"{o}lder inequality last,  then
\begin{eqnarray}\label{3.25}
\mE\sup_{0\leq t\leq T}|Y_t^1-Y_t^2|^2\leq  C\int_0^T(1+\kappa(t)) \mE|Y_t^1-Y_t^2|^2dt +\Big[\int_0^T \kappa(t)\mE|Y_t^1-Y_t^2|^4dt\Big]^{\frac12}=0,
\end{eqnarray}
where in the last identity in (\ref{3.25}) we have used (\ref{3.24}). Therefore, pathwise uniqueness of solutions holds for (\ref{3.22}) and thus holds for (\ref{1.1}).

It remains to show the H\"{o}lder continuity. By virtue of Lemma \ref{lem3.2}, the scaling
transformation and the continuity of $X$ in $t$, it suffices  to check that for every $p\geq2$, $x,y\in\mathbb{R}^d$,
\begin{eqnarray}\label{3.26}
\mathbb{E}[\sup_{0\leq t\leq 1}|Y_t(x)-Y_t(y)|^p]\leq C|x-y|^p.
\end{eqnarray}

Let $\kappa$ be given in (\ref{3.23}). From (\ref{3.22}),
by employing the It\^{o} formula, then
\begin{eqnarray}\label{3.27}
\mathbb{E}|Y_t(x)-Y_t(y)|^p\leq |x-y|^p+C\int_0^t|Y_r(x)-Y_r(y)|^pdr +\int_0^t\kappa(r)|Y_r(x)-Y_r(y)|^pdr.
\end{eqnarray}
Applying  the Gr\"{o}nwall inequality, for every $p\geq2$, we get
\begin{eqnarray}\label{3.28}
\sup_{0\leq t\leq 1}\mathbb{E}|Y_t(x)-Y_t(y)|^p\leq C|x-y|^p.
\end{eqnarray}

By (\ref{3.22}), the Doob and BDG inequalities, then
\begin{eqnarray}\label{3.29}
&&\mathbb{E}\sup_{0\leq t \leq 1}|Y_t(x)-Y_t(y)|^p\nonumber\\&\leq & |x-y|^p+C\mathbb{E}\int_0^1(1+\kappa(r))|Y_r(x)-Y_r(y)|^pdr
\nonumber\\&& +C\Big[\mathbb{E}\Big|\int_0^1|Y_r(x)-Y_r(y)|^{2p-2}
|\tilde{\sigma}(r,Y_r(x))-\tilde{\sigma}(r,Y_r(y))|^2dr\Big]^{\frac{1}{2}}
\nonumber\\&\leq& |x-y|^p+C\mathbb{E}\int_0^1(1+\kappa(r))|Y_r(x)-Y_r(y)|^pdr
+C\Big[\mathbb{E}\int_0^1\kappa(r)|Y_r(x)-Y_r(y)|^{2p}dr\Big]^{\frac{1}{2}}.
\end{eqnarray}
Observing that (\ref{3.28}) holds for every $p\geq2$, from (\ref{3.29}),  (\ref{3.26}) holds.

$(iii)$ By the relationship between $X_t$ and $Y_t$, to prove the homeomorphism property of $X_t(\cdot)$ it needs to prove that for almost all $\omega\in\Omega$, every $t>0$, $y\rightarrow Y_t(y)$ is a homeomorphism on $\mathbb{R}^d$. Due to \cite[Theorem 4.5.1]{Kun} and the fact
\begin{eqnarray*}
\sup_{y\in\mathbb{R}^d}\mathbb{E}|Y_t(y)-Y_r(y)|^p\leq C|t-r|^{\frac{p}{2}},  \quad 0\leq t,r\leq T,
\end{eqnarray*}
we should prove that: for every $\tau\in \mathbb{R}$ and $x,y\in \mathbb{R}^d$ ($x\neq y$)
\begin{eqnarray}\label{3.30}
\mathbb{E}\sup_{0\leq t\leq T}(1+|Y_t(y)|^2)^\tau\leq C(1+|y|^2)^\tau
\end{eqnarray}
and
\begin{eqnarray}\label{3.31}
\sup_{0\leq t\leq T}\mathbb{E}|Y_t(x)-Y_t(y)|^{2\tau}\leq C|x-y|^{2\tau}.
\end{eqnarray}
Since $\tilde{b}$ and $\tilde{\sigma}$ are bounded, (\ref{3.30}) is obvious.
It remains to calculate (\ref{3.31}). For $\epsilon>0$, if one
chooses $F(x)=f^\tau(x)=(\epsilon+|x|^2)^\tau$ and set
$Y_t(x,y):=Y_t(x)-Y_t(y)$, then by utilising the It\^{o} formula,
\begin{eqnarray}\label{3.32}
F(Y_t(x,y))&=&F(Y_0(x,y))+2\tau\int_0^tf^{\tau-1}(Y_r(x,y))\langle Y_r(x,y), \tilde{b}(r,Y_r(x))-\tilde{b}(r,Y_r(y))\rangle dr\nonumber\\&&+
2\tau\int_0^tf^{\tau-1}(Y_r(x,y))\langle Y_r(x,y), (\tilde{\sigma}(r,Y_r(x))-\tilde{\sigma}(r,Y_r(y)))dW_r\rangle
\nonumber\\&&+\tau\sum_{i,j=1}^d\int_0^tf^{\tau-2}(Y_r(x,y))[f(Y_r(x,y))\delta_{i,j}+
2(\tau-1)Y_{i,r}(x,y)Y_{j,r}(x,y)]
\nonumber\\&&\quad\times[\tilde{\sigma}_{i,k}(r,Y_r(x))-
\tilde{\sigma}_{i,k}(r,Y_r(y))][\tilde{\sigma}_{j,k}(r,Y_r(x))
-\tilde{\sigma}_{j,k}(r,Y_r(y))]dr
\nonumber\\&\leq& F(Y_0(x,y))+C|\tau|\int_0^tF(Y_r(x,y))dr+C|\tau(\tau-1)|\int_0^t\kappa(r)
F(Y_r(x,y))dr\nonumber\\&&+
2\tau\int_0^tf^{\tau-1}(Y_r(x,y))\langle Y_r(x,y), (\tilde{\sigma}(r,Y_r(x))-\tilde{\sigma}(r,Y_r(y)))dW_r\rangle,
\end{eqnarray}
where $\kappa$ is given in (\ref{3.23}), and
\begin{eqnarray*}
\delta_{i,j}=\left\{\begin{array}{ll}
1, \ \ \mbox{when} \ \ i=j, \\
0, \  \ \mbox{when} \ \ i\neq j.
\end{array}\right.
\end{eqnarray*}
Thanks to (\ref{3.32}) and the Gr\"{o}nwall inequality, one arrives at
\begin{eqnarray*}
\sup_{0\leq t\leq T}\mathbb{E}[\epsilon+|Y_t(x)-Y_t(y)|^2]^\tau\leq C[\epsilon+|x-y|^2]^\tau.
\end{eqnarray*}
By letting $\epsilon\downarrow 0$, then (\ref{3.31}) holds.

Similarly, to prove the weak differentiability of $X_t(\cdot)$, it only needs to show
\begin{eqnarray}\label{3.33}
\lim_{\delta\rightarrow 0}\frac{Y_\cdot(y+\delta e_i)-Y_\cdot(y)}{\delta}
\end{eqnarray}
exists in $L^2(\Omega\times (0,T))$. Set
$Y_t^\delta(y):=Y_t(y+\delta e_i)-Y_t(y)$, then by (\ref{3.22})
\begin{eqnarray}\label{3.34}
Y_t^\delta(y)&=&\delta e_i+\int_0^t[\tilde{b}(r,Y_r(y+\delta e_i))-\tilde{b}(r,Y_r(y))]dr+\int_0^t[\tilde{\sigma}(r,Y_r(y+\delta e_i))-\tilde{\sigma}(r,Y_r(y))]dW_r\nonumber\\&=&\delta e_i+\int_0^1\int_0^t\nabla\tilde{b}
(r,sY_r(y+\delta e_i)+(1-s)Y_r(y))Y_r^\delta(y)drds
\nonumber\\&&+\int_0^1\int_0^t[\tilde{\sigma}(r,sY_r(y+\delta e_i)+(1-s)Y_r(y))]Y_r^\delta(y)dW_rds.
\end{eqnarray}

By virtue of BDG's inequality, we achieve from (\ref{3.34}) that
\begin{eqnarray*}
\mathbb{E}|Y_t^\delta(y)|^2\leq 2|\delta|^2 +C\mathbb{E}\int_0^t|Y_r^\delta(y)|^2dr+\mathbb{E}
\int_0^t\kappa(r)|Y_r^\delta(y)|^2dr,
\end{eqnarray*}
which suggests that
\begin{eqnarray*}
\mathbb{E}\int_0^T\Big|\frac{|Y_t^\delta(y)|}{\delta}\Big|^2dt\leq C.
\end{eqnarray*}
Notice that $L^2((0,T)\times\Omega)$ is reflexive, there is a subsequence still denoted by itself such that $Y_\cdot^\delta$ converges weakly to an element $\tilde{Y}$ in $L^2((0,T)\times\Omega)$.  Then by applying the resonance theorem,
\begin{eqnarray*}
\mathbb{E}\int_0^T|\tilde{Y}|^2dt\leq \liminf_{\delta\rightarrow 0} \mathbb{E}\int_0^T\Big|\frac{|Y_t^\delta(y)|}{\delta}\Big|^2dt\leq C.
\end{eqnarray*}
The desired result follows. $\Box$

\begin{remark} \label{rem3.1} When $d=1$, we also derive the pathwise uniqueness without assuming the
Sobolev differentiability on $\sigma$ if $\alpha\geq1/2$. In this case, $\tilde{b}\in \mathcal{C}([0,T];\mathcal{C}_b^1(\mathbb{R})),
\tilde{\sigma}\in L^\infty(0,T;{\mathcal C}^\alpha(\mathbb{R}))$,
by Yamada-Watanabe's theorem (see \cite{YW}), the pathwise uniqueness
holds for SDE (\ref{3.22}) and thus for (\ref{1.1}).
\end{remark}

\section{Conclusions}\label{sec4}\setcounter{equation}{0}

In recent years, people have made broad research about the uniqueness of strong solutions
for the stochastic differential equation
\begin{eqnarray}\label{4.1}
dX_t=b(t,X_t)dt+\sigma(t,X_t)dW_t, \ t\in (0,T], \
X_0=x\in{\mathbb R}^d,
\end{eqnarray}
with non-Lipschitz coefficients. However, most of these works are concentrated on the drift which is integrable in time with the integrable index $q\in (2,+\infty]$. There are few research works concerned with $q=2$. In this study, we have established the existence, uniqueness H\"{o}lder continuity and weak differentiability of strong solutions only assuming $b\in L^2(0,T;\cC^\alpha_b(\mR^d;\mR^d))$. Compared with the existing research, these results are new.

To prove the results, we use the It\^{o}-Tanaka trick  to
transform the SDE (\ref{4.1}) with singular coefficients to an equivalent new SDE with Lipschitz coefficients via a non-singular diffeomorphism $\Phi(t,x)=x+U(t,x)$, where $U(T-t,x)=:V(t,x)$ satisfies a vector-valued parabolic partial differential equation of second order which has the form:
\begin{eqnarray}\label{4.2}
\left\{\begin{array}{ll}
\partial_{t}V(t,x)=\frac{1}{2}\sum_{i,j=1}^da_{i,j}(t,x)\partial^2_{x_i,x_j} V(t,x)+b(t,x)\cdot \nabla V(t,x)
\\ \qquad\qquad\qquad+b(t,x)-\lambda V(t,x), \ (t,x)\in (0,T)\times {\mathbb R}^d, \\
V(0,x)=0, \  x\in{\mathbb R}^d.  \end{array}\right.
\end{eqnarray}
To accomplish the goal, there are two things we need to do. The first one is choosing a proper function space on which the It\^{o} formula is applicable, and the second one is the boundedness estimate for the gradient of $U$. We accomplish the first issue by fetching $L^2(0,T;\cC^{2,\alpha}_b(\mR^d))\cap W^{1,2}(0,T;\cC^{\alpha}_b(\mR^d))$ as the workspace. When $b\in L^\infty(0,T;\mathcal{C}^{\alpha}_b(\mathbb{R}^d;\mathbb{R}^d))$,
these estimates for solutions have been established by Krylov
\cite{Kry02}.
Noticing that, here we only assume that
$b\in L^2(0,T;\mathcal{C}_b^\alpha(\mathbb{R}^d;\mathbb{R}^d))$,
so we should extend Krylov's result. We found these estimates for solutions of $V$ by using the fundamental solution, and the key point is to use the heat kernel estimates for the following scalar parabolic partial differential equation:
\begin{eqnarray}\label{4.3}
\left\{\begin{array}{ll}
\partial_{t}u(t,x)=\frac{1}{2}\sum_{i,j=1}^da_{i,j}(t,x)\partial^2_{x_i,x_j} u(t,x)+b(t,x)\cdot \nabla u(t,x), \ (t,x)\in (0,T)\times {\mathbb R}^d, \\
u(0,x)=0, \  x\in{\mathbb R}^d.  \end{array}\right.
\end{eqnarray}
When $a_{i,j}$ are in $L^\infty(0,T;\mathcal{C}^{\alpha}_b(\mathbb{R}^d))$ and $b$ is in $L^q(0,T;\cC^\alpha_b(\mR^d;\mR^d))$ $(q>2)$, the existence, two-sides estimates and gradient estimates for heat kernel have been established by Chen, Hu, Zhang and Xie \cite[Theorem 2.3]{Chen17} since now $b$ is in the generalized Kato class $\mK_2$ (see \cite[Definition 2.5]{Chen17}). Unfortunately, when $b$ is only square integrable in time, for a given time $t_0\in (0,T)$, if one chooses $b(t,x)=(t_0-t)^{-1/2}1_{0\leq t< t_0}|\log(t_0-t)|^{-1} (x\in\mR)$ it is not in $\mK_2$. New ideas are needed to approach this problem. In this paper, we use the heat kernel for (\ref{4.3}) with $b=0$ and combine it with Banach's fixed point theorem to prove the Schauder estimates of solutions for (\ref{4.3}) with a nonhomogeneous term. Then using the same techniques to establish the gradient estimate for $u$. These results are new as well.

\bigskip\noindent
\textbf{\large{Acknowledgements.}} The first author was partially supported by the National Natural Science Foundation of China (11901442). The second author was partially supported by the research funding project of Guizhou Minzu University (GZMU[2019]QN04). The third author was partially supported by the National Natural Science Foundation of China (11501577).

\end{document}